\documentclass[a4paper]{amsart}

\usepackage{amsfonts,amsthm,amssymb,amsmath,amscd}

\newtheorem{lem}{Lemma}[section]
\newtheorem{prop}{Proposition}[section] 
\newtheorem{cor}{Corollary}[section] 
\newtheorem{rem}{Remark}[section] 
\newtheorem{tm}{Theorem}[section] 

\newcommand{\ZZ}{\mathbb{Z}}
\newcommand{\QQ}{\mathbb{Q}}

\newcommand{\CC}{\mathbb{C}}
\newcommand{\PP}{\mathbb{P}}

\begin{document}

\title{Algebraicity of some Weil Hodge classes}            

\author{Kenji Koike}

\maketitle

\begin{abstract} 
We show that the Prym map for 4-th cyclic \'etale covers of 
curves of genus 4 is a dominant morphism to a Shimura variety for a family 
of Abelian 6-folds of Weil type. According to the result of Schoen, 
this implies algebraicity of Weil classes for this family.
\end{abstract}

\section{Introduction}
The Hodge conjecture is still open even for Abelian varieties. 
Because the conjecture is true for all projective 3-folds, the first targets 
are 4-folds. In this dimension, 
we have non-trivial examples of Hodge classes for special Abelian varieties 
that are called of Weil type.
\par
In \cite{Sc1} and \cite{Sc2}, Schoen constructed algebraic cycles on 
generalized Prym varieties for cyclic covers that give Weil classes, 
and he proved algebraicity of Weil classes for a family of Abelian 
6-folds of Weil type for $k = \QQ(\sqrt{-3})$ with some $\delta$ 
by showing the denseness of Prym varieties in the family. 
\par
This method works also for Weil 4-folds and 6-folds for $\QQ(\sqrt{-1})$ 
with $\delta = 1$. For the 4-dimensional case, van Geemen gave another proof 
in \cite{vG2}. 
\par
In this note, we consider the 6-dimensional case. 
The problem is to show the dominantness of the associated Prym map. 
We construct genus 13 curves $C_{13}$ in $(\PP^1)^4$ that are invariant 
under a cyclic permutation $\sigma$ of factors of $(\PP^1)^4$. So $\sigma$ 
acts on $C_{13}$, and we show that this action is fixed point free. For 
the covering $C_{13} \rightarrow C_{13} / \left< \sigma \right >$, we compute 
the codifferential map of the Prym map explicitly. 
\vskip3mm
{\bf acknowledgements}
This work was done at Johann Wolfgang Goethe-Universit\"at while the 
author was supported by Alxaander von Humboldt Stiftung. He is grateful to 
Prof. J. Wolfart for many suggestions and his hospitality. He also 
thanks Prof. B. van Geemen at Universita Milano for explaining some 
detailed arguments.
\section{The Prym construction of Abelian varieties of Weil type}
In this section, we explain our problem and state the main theorem.
\par
Let us recall the definition of Abelian varieties of Weil type 
(see \cite{vG1} and \cite{We}). 
Let $A$ be an Abelian $2n$-fold ($n \geq 2$) with a polarization $E$ 
and let $\phi : k \rightarrow \mathrm{End}^0(A)$ be an inclusion 
of an imaginary quadratic field $k  = \QQ(\sqrt{-d})$. We assume that 
\[
E(\phi(\sqrt{-d})_* x, \phi(\sqrt{-d})_* y) = d E(x,y).
\] 
Then we say that $A$ is of Weil type, if the multiplicities of 
eigenvalues $\sqrt{-d}$ and $-\sqrt{-d}$ of the action of 
$\phi (\sqrt{-d})_*$ on the tangent space $T_0A$ are equal to $n$. 
Then $H_1(A, \QQ)$ has a structure of $k$-module and
\[ 
H(x,y) = E(x,\phi(\sqrt{-d})_* y) + \sqrt{-d} E(x,y)
\]
gives a Hermitian form on $H_1(A, \QQ)$ of the signature $(n,n)$. 
The class $\delta = \det H$ mod $(k^*)^2$ gives an isogenus invariant, and 
we call this the discriminant of $(A, E, \phi)$.
For these Abelian $2n$-folds, Weil constructed a subspace
\[
W(A) = \wedge^n_k H^1(A, \QQ) 
\subset H^{n,n}(A) \cap H^{2n}(A, \QQ) 
\]
of the Hodge classes ($\dim_{\QQ} W(A) = 2$), 
and elements in $W(A)$ are called Weil classes. The special Munford-Tate 
group for a general $A$ is a special unitary group of the signature $(n, n)$, 
and in this case we have
\[
H^{n,n}(A) \cap H^{2n}(A, \QQ) = D^n(A) \oplus W(A),
\]
where $D^{\bullet}(A)$ is the subspace of classes generated by divisors. 
To prove the Hodge conjecture, therefore, 
we need algebraic cycles that never come from divisors.
\par
In special cases, the setting up of the problem is established in \cite{Sc1}
(see also \cite{vG2}).
We consider a curve $C_{13}$ of genus $13$ which is a 4-th cyclic 
\'etale cover of a curve $C_4$ of genus $4$. Then we have \'etale double 
cover $C_{13} \rightarrow C_7$ of a curve $C_7$ of genus $7$ as 
the intermediate cover. Let us consider the Prym variety 
$P = Prym(C_{13}/C_7)$
that is the connected component of the kernel of the norm map 
$Nm : J(C_{13}) \rightarrow J(C_7)$ including $0 \in J(C_{13})$. 
This is a principally polarized 6-dimensional Abelian variety.
The Galois group $Gal(C_{13}/C_4) \cong \ZZ / 4 \ZZ$ acts on $P$, 
and $P$ becomes Abelian variety of Weil type for $\QQ(\sqrt{-1})$ by 
this action. It is known that the discriminant $\delta$ of $P$ is $1$ 
(see \cite{vG2}), and the Weil classes $W(P)$ are generated by algebraic 
cycles (see \cite{Sc1}).
\par
Let $\mathcal{M}_{13/4}$ be the moduli space of 4-th cyclic covers 
$C_{13} \rightarrow C_4$ and $\mathcal{A}_6$ be the moduli space of 
principally polarized Abelian varieties of dimension $6$. Then we have 
the Prym map
\[
Pr : \mathcal{M}_{13/4} \longrightarrow \mathcal{H}_6 \subset \mathcal{A}_6,
\quad \{C_{13} \rightarrow C_7 \rightarrow C_4\} \mapsto Prym(C_{13} / C_7),
\]
where $\mathcal{H}_6$ is the Shimura variety of dimension $9$ 
given in \cite{vG2}. Therefore the Weil classes $W(A)$ is generated by 
algebraic cycles for a general $A \in \mathcal{H}_6$ if the image of $Pr$ is 
Zariski dense in $\mathcal{H}_6$. Because the moduli space 
$\mathcal{M}_{13/4}$ is finite over the moduli space $\mathcal{M}_4$ of 
genus 4 curves and $\mathrm{dim} \mathcal{M}_4 = 9$, 
it is enough to prove the dominantness of 
$Pr$ if we show that the differential (equivalently, the codifferential) of 
$Pr$ at some point of $\mathcal{M}_{13/4}$ is an isomorphism.
\par
Let us take $\Pi = \{X \rightarrow Y \rightarrow Z \} 
\in \mathcal{M}_{13/4}$ and the line bundle $L$ on $Z$ which gives 
$\pi:X \rightarrow Z$. Then we have the decomposition
\[
\pi_* \omega_X = \omega_Z \oplus (\omega_Z \otimes L) \oplus 
(\omega_Z \otimes L^2) \oplus (\omega_Z \otimes L^3).
\]
Replacing moduli spaces by them with a level structure if necessary, 
cotangent spaces at $\Pi$ and at the intermediate Prym variety $P$ of 
$\Pi$ are represented by
\begin{align*}
&T^*_P \mathcal{H}_6 = (T^*_P \mathcal{A}_6)^{Gal(X/Z)} =
H^0(Z, \omega_Z \otimes L) \otimes H^0(Z, \omega_Z \otimes L^3), \\
&T^*_{\Pi} \mathcal{M}_{13/4} = T^*_{Z} \mathcal{M}_4 = 
H^0(Z, \omega_Z \otimes \omega_Z),
\end{align*}
and the codifferential map of $Pr$ at $\Pi$ is given the multiplication map
\begin{align} \label{mu}
\mu : H^0(Z, \omega_Z \otimes L) \otimes
H^0(Z, \omega_Z \otimes L^3) \longrightarrow 
H^0(Z, \omega_Z \otimes \omega_Z).
\end{align}
\par
Now we state our result:
\begin{tm}
There exist a $4$-th cyclic \'etale cover $\pi : X \rightarrow Z$ of genus 
$4$ curve $Z$ such that the multiplication map $\mu$ in (\ref{mu}) is an 
isomorphism, and therefore the Prym map $Pr$ is dominant.
\end{tm} 
By the specialization argument in \cite{Sc1}, we know that
\begin{cor}
The Weil classes are generated by algebraic cycles for Abelian $6$-folds of 
Weil type for $\QQ(\sqrt{-1})$ with $\delta = 1$.
\end{cor}
\begin{rem}
By the Proposition 10 in \cite{Sc2}, we see that the Weil classes are 
generated by algebraic cycles for all Abelian $4$-folds of Weil type 
for $\QQ(\sqrt{-1})$.
\end{rem}
\section{Complete intersections in $(\PP^1)^4$ and  a cyclic permutation}
In this section we construct a genus $13$ curve $C_{13}$ in $(\PP^1)^4$ 
with  a fixed point free automorphism $\sigma$ of order $4$, and we show that 
the natural projection $C_{13} \rightarrow C_{13}/ \left< \sigma \right>$ 
satisfies the required condition.
\par
Let $[s_0:s_1] \times [t_0:t_1] \times [x_0:x_1] \times [y_0:y_1]$ be the 
coordinate of $(\PP^1)^4$. For the simplicity, we denote this by 
$(s,t,x,y)$ with $s = s_1/s_0$ and so on. Let $\sigma$ be a cyclic 
permutation $(s,t,x,y) \mapsto (y,s,t,x)$ on $(\PP^1)^4$. 
Then $\sigma$ acts on the vector space 
\[
V = H^0((\PP^1)^4, \otimes_{i = 1}^4 p_i^*\mathcal{O}_{\PP^1}(1)) = 
\oplus \CC s_i t_j x_k y_l, \quad i,j,k,l \in \{0,1\},
\]
where $p_i : (\PP^1)^4 \rightarrow \PP^1$ is the $i$-th projection. 
Let $V(\alpha)$ be the eigenspace for the eigenvalue $\alpha$. 
We have the following basis of $V(\alpha)$:
\begin{equation} \label{base}
\begin{split}
V(1): \ &a_1 = s + t + x + y, \quad a_2 = st + tx + xy + ys, \\
&a_3 = txy + sxy + sty + stx, \quad 
a_4 = sx + ty, \\ 
&a_5 = stxy, \quad a_6 = 1,\\
V(-1): \ &b_1 = s - t + x - y, \quad b_2 = st - tx + xy - ys, \\
&b_3 = txy - sxy + sty - stx, \quad b_4 = sx - ty, \\
V(i): \ &c_1 = s -i t - x +i y, \quad c_2 = st -i tx - xy +i ys, \\
&c_3 = txy -i sxy - sty +i stx, \\
V(-i): \ &d_1 = s +i t - x -i y, \quad d_2 = st +i tx - xy -i ys, \\
&d_3 = txy +i sxy - sty -i stx,
\end{split}
\end{equation}
where we identified a multi-homogeneous polynomial $f \in V$ with 
the polynomial $f/(s_0 t_0 x_0 y_0)$ in affine coordinates. 
Note that every element in $V(1)$ is invariant also for the involution 
\[
\tau : (\PP^1)^4 \longrightarrow (\PP^1)^4, \quad (s,t,x,y) 
\mapsto (x,t,s,y),
\]
namely, they are invariant under the dihedral group 
$G = \left< \sigma, \tau \right>$.
\begin{lem} \label{lem1}
The linear system $V(1)$ gives a morphism 
$\varphi : (\PP^1)^4 \rightarrow F \subset \PP^5$ of generic degree $8$, 
where $F$ is the cubic hypersurface
\begin{align}
a_1^2 a_4 - a_1 a_3 a_4 + a_2 a_4^2 - 4 a_2 a_5 a_6 + a_3^2 a_6 = 0.
\end{align}
So the fiber $\varphi^{-1}(P)$ of a generic point $P \in F$ is a $G$-orbit.
\end{lem}
\begin{proof}
The morphism
\[
(\PP^1)^4 \longrightarrow \PP^4, \quad (s,t,x,y) \mapsto
[a_1:a_2 + a_4:a_3:a_5:a_6]
\]
induces an isomorphism $(\PP^1)^4/S_4 \rightarrow \PP^4$ because they 
are the fundamental symmetric polynomials. Hence $V(1)$ has no base points, 
and $\varphi$ is a finite morphism onto the image. 
We can check that  $a_i$'s satisfy the above cubic equation, and 
that the morphism $F \rightarrow \PP^4$ is of generic degree $3$. 
Now we see that $\varphi$ is of degree $8$ since $|S_4| = 24$. 
\end{proof}
\begin{lem} \label{lem2}
Let $C$ be a hyperelliptic curve of genus $4$, $L$ be a non-trivial 
line bundle such that $L^2 = \mathcal{O}_C$. Then the associated map 
$\varphi_{\omega_C \otimes L} : C \rightarrow \PP^2$ 
defined by $H^0(C, \ \omega_C \otimes L)$ satisfies one of 
the following conditions.
\par
(1). $\varphi_{\omega_C \otimes L}$ is a rational map to a conic in $\PP^2$.
\par
(2). $\varphi_{\omega_C \otimes L}$ is a birational map.
\end{lem}
\begin{proof}
Because $L$ is isomorphic to $\mathcal{O}_C(P_1 - P_2)$ or to 
$\mathcal{O}_C(P_1 + P_2 - P_3 - P_4)$ 
where $\pi(P_i)$'s are distinct blanch points of the double cover 
$\pi:C \rightarrow \PP^1$, it is enough if we consider these two cases. 
\par 
Let $y^2 = f(x)$ be the equation of $C$  and $\pi(P_i) = \lambda_i$. 
Then a basis of $H^0(C, \ \omega_C \otimes L)$ is given by 
\[
(x-\lambda_2) \frac{dx}{y}, \quad (x-\lambda_2)^2 \frac{dx}{y}, \quad 
(x-\lambda_2)^3 \frac{dx}{y}, \quad
\]
for $L = \mathcal{O}_C(P_1 - P_2)$, and  this is the case (1) in the 
assertion. In the case of $L = \mathcal{O}_C(P_1 + P_2 - P_3 - P_4)$, we 
can take the following basis of $H^0(C, \ \omega_C \otimes L)$
\[
\eta_1 = (x-\lambda_3)(x-\lambda_4) \frac{dx}{y}, \quad 
\eta_2 = (x-\lambda_3)(x-\lambda_4)^2 \frac{dx}{y}, \quad 
\eta_3 = \frac{dx}{(x-\lambda_1)(x-\lambda_2)}.
\]
Because $\eta_2 / \eta_1 = x - \lambda_4$ and $\eta_3 / \eta_1 = g(x) y$ with 
a rational function $g(x)$, we see that $\varphi_{\omega_C \otimes L}$ 
is birational in this case.
\end{proof}
\begin{prop} \label{prop1}
For general elements $f_1, f_2, f_3 \in V(1)$,
the complete intersection $f_1 = f_2 = f_3 = 0$ defines 
a smooth curve $X$ of genus $13$, and we have an isomorphism of 
vector spaces
\[
H^0(X, \omega_X) \cong 
V(1)/(\CC f_1 \oplus \CC f_2 \oplus \CC f_3) \oplus V(-1) \oplus V(i) 
\oplus V(-i). 
\]
The cyclic permutation $\sigma$ acts on $X$ without fixed point.
Hence we have a \'etale cyclic cover 
$X \rightarrow Y = X/ \left< \sigma^2 \right> \rightarrow 
Z = X/ \left< \sigma \right>$,
and an isomorphism $H^0(Z, \omega_Z) \cong V(-1)$. The genus $4$ curve 
$Z$ is not hyperelliptic.
\end{prop}
\begin{proof}
Because $V(1)$ is base point free and the divisor given by a general 
$f \in V(1)$ is reduced, the curve 
$X = \{f_1 = f_2 = f_3 = 0\}$ is smooth for general $f_1, f_2, f_3 \in V(1)$. 
By the Adjunction formula, we know that the restriction of $V$ on $X$ gives 
the canonical class and the genus of $X$ is $13$. 
\par
Obviously $\sigma$ acts on $X$, and the fixed points of 
$\sigma^2$ on $(\PP^1)^4$ is 
\[
\Delta = \{ (s,t,s,t) \in (\PP^1)^4 \} \cong \PP^1 \times \PP^1. 
\]
The restriction of the basis of $V(1)$ in (\ref{base}) to 
$\Delta$ is given by
\begin{equation}
\begin{split}
&a_1 = s_0t_0(s_1 t_0 + s_0 t_1), \quad 
a_2 = s_0t_0s_1t_1, \quad a_3 = s_1t_1(s_1t_0 + s_0t_1), \\ 
&a_4 = s_1^2t_0^2 + s_0^2t_1^2, \quad a_5 = s_1^2t_1^2, 
\quad a_6 = s_0^2t_0^2,
\end{split}
\end{equation}
up to constant, with coordinates $[s_0:s_1] \times [t_0:t_1]$. They have 
no base point and we see that $C \cap \Delta = \phi$ for general $f_i$'s. 
\par
Because $V(-1)$ is the unique $4$-dimensional eigenspace for the action of 
$\sigma$, we may identify $V(-1)$ with $H^0(Z, \omega_Z)$. 
\par
Now let $L$ be the line bundle corresponding to $\pi: Y \rightarrow Z$. 
We see that 
\[
H^0(Z, \omega_Z \otimes L) \cong 
V(1)/(\CC f_1 \oplus \CC f_2 \oplus \CC f_3).
\]
By Lemma (\ref{lem1}), the image of $\varphi_{\omega_Z \otimes L}$ 
is a cubic curve $E$ which is a section of $F$ by $\PP^3 \subset \PP^5$, 
and this is not any case in Lemma (\ref{lem2}). 
Therefore $Z$ is not hyperelliptic.
\end{proof}
\begin{rem}
We can check that the singular locus of $F$ is $1$-dimensional, and that 
the section of $F$ by a generic $\PP^3$ in $\PP^5$ is a smooth cubic 
curve $E$. The Prym canonical map $\varphi_{\omega_Z \otimes L}$ is just the 
natural map $Z \rightarrow E = Z / \left< \tau \right>$, so our curve $Z$ 
is bi-elliptic.
\end{rem}
\par
Let $f_1, f_2, f_3$ and $X$ be as in Proposition (\ref{prop1}), 
$\pi : X \rightarrow Z = X/ \left< \sigma \right>$ be the quotient map 
and $L$ be the line bundle on $Z$ corresponding to $\pi$. Because $Z$ 
is not hyperelliptic, the multiplication map
\[
Sym^2 H^0(Z, \omega_Z) \longrightarrow H^0(Z, \omega_Z \otimes \omega_Z)
\]
is a surjection by Max Noether's Theorem (see \cite{ACGH}). By Riemann-Roch 
Theorem, we see that this map has the $1$-dimensional kernel. We denote 
a generator of this kernel by $Q$, and we use the same symbol for the 
corresponding element in $Sym^2V(-1)$. Namely, we have isomorphisms
\begin{align*}
Sym^2V(-1) / \CC Q \cong H^0(Z, \omega_Z \otimes \omega_Z), \\
V(i) \otimes V(-i) \cong H^0(Z, \omega_Z \otimes L) \otimes 
H^0(Z, \omega_Z \otimes L^3).
\end{align*}
Therefore the map $\mu$ in (\ref{mu}) defines the induced multiplication map
\begin{align}
m : V(i) \otimes V(-i) \longrightarrow Sym^2V(-1) / \CC Q, \quad 
c \otimes d \mapsto cd \ \text{mod} \ (f_1, f_2, f_3) 
\end{align}
by the above identification.
(The map $m$ is well-defined only modulo $f_1 = f_2 = f_3 = 0$.) 
Now the bijectivity 
of $\mu$ is equivalent to the linear independence of 
$\{c_i d_j \}_{1 \leq i,j \leq 3}$ modulo $f_1 = f_2 = f_3 = 0$. 
\par
Let us show the linear independence. First of all, we have 
the following quadric equations (these are 
a part of the Segre relations, and we can check them with a computer)
\begin{equation} \label{eqb}
\begin{split}
&b_1^2 = a_1^2 - 4 a_2 a_6, \quad 
b_2^2 = a_2^2 - 4 (a_1 a_3 - a_2 a_4) + 16 a_5 a_6, \\ 
&b_3^2 = a_3^2 - 4 a_2 a_5, \quad
b_1 b_3 = a_1 a_3 - 2 a_2 a_4, \quad 
b_1 b_4 = a_1 a_4 - 2 a_3 a_6, \\
&b_3 b_4 = - a_3 a_4 + 2 a_1 a_5, \quad b_4^2 = a_4^2 - 4 a_5 a_6,
\end{split}
\end{equation}
and
\begin{equation} \label{eqcd}
\begin{split}
&c_1 d_1 = a_1^2 - 2 a_2 a_6 - 4 a_4 a_6, \quad
c_3 d_3 = a_3^2 - 2 a_2 a_5 - 4 a_4 a_5, \\
&c_2 d_2 = a_2^2 - 2 a_1 a_3 + 2 a_2 a_4, \\
&\frac{1}{2}(c_1 d_3 + c_3 d_1) = - a_1 a_3 + a_2 a_4  + 8 a_5 a_6, \quad
\frac{i}{2} (c_1 d_3 - c_3 d_1) = b_2 b_4, \\
&\frac{1}{1-i}(c_1 d_2 -i c_2 d_1) = a_1 a_2 - 4 a_3 a_6, \quad
\frac{1}{1+i}(c_1 d_2 +i c_2 d_1) = b_1 b_2, \\
&\frac{1}{1+i}(c_2 d_3 + i c_3 d_2) = -a_2 a_3 + 4 a_1 a_5, \quad
\frac{-i}{1+i}(c_2 d_3 - i c_3 d_2) = b_2 b_3.
\end{split}
\end{equation}
Without a loss of generality, we may assume that our equations 
$f_1 = 0$, $f_2 = 0$ and $f_3 = 0$ are given by 
\begin{equation} \label{eqa}
\begin{split}
&a_4 = A_1 a_1 + A_2 a_2 + A_3 a_3, \quad
a_5 = B_1 a_1 + B_2 a_2 + B_3 a_3, \\
&a_6 = C_1 a_1 + C_2 a_2 + C_3 a_3,
\end{split}
\end{equation}
with coefficients $A_i, B_i, C_i \in \CC$. Substituting them, we can 
eliminate $a_4, a_5, a_6$ in equations in (\ref{eqb}). Then 
each product $b_i b_j$ in (\ref{eqb}) is a linear combination of six 
elements
\begin{align} \label{6a}
a_1^2, \quad a_2^2, \quad a_3^2, \quad a_1 a_2, \quad a_1 a_3, 
\quad a_2 a_3.
\end{align}
Since there are seven elements in (\ref{eqb}), we have a non-trivial linear 
relation of $b_i b_j$'s in (\ref{eqb}) which gives the unique vanishing 
quadric $Q$. Therefore $a_i a_j$'s in (\ref{6a}) and 
$b_1 b_2, b_2 b_3, b_2 b_4$ give a base of the vector space 
$Sym^2V(-1)/ \CC Q $ (modulo equations in (\ref{eqa})). 
\par
Now let us consider a base of the vector space $V(i) \otimes V(-i)$ given 
in (\ref{eqcd}). Eliminating $a_4, a_5, a_6$, we obtain a vector equation
$\gamma = M \alpha$, where 
\begin{align*}
\gamma = {}^t(c_1 d_1,\quad c_2 d_2, \quad c_3 d_3, \quad 
\frac{1}{1-i}(c_1 d_2 - ic_2 d_1), \quad \frac{1}{2}(c_1 d_3 + c_3 d_1), \\
\quad \frac{1}{1+i}(c_2 d_3 + ic_3 d_2)), \\
\alpha = {}^t(a_1^2, \quad a_2^2, \quad a_3^2, \quad a_1 a_2, 
\quad a_1 a_3, \quad a_2 a_3),
\end{align*}
and $M$ is a matrix of polynomials in $A_i, B_i, C_i$. If we have 
$\det M \ne 0$, we can conclude that $c_i d_j$'s form a base of 
$Sym^2V(-1)/ \CC Q $, and we finish the proof of theorem. 
\par
Let us consider the case that all $A_i, B_i, C_i$ are $0$. 
Then $a_4 = a_5 = a_6 = 0$, and obviously we have $\det M = 1$. 
Therefore $\det M$ is not identically zero as a polynomial of 
$A_i, B_i, C_i$.

\end{document}